\newcommand{\elsevierformat}{
  \documentclass{elsart}
  \newenvironment{proof}{\begin{pf}}{\end{pf}}
  \newcommand{\bstyle}{elsart-num}
  \newcommand{\maketitle}{}

  \newenvironment{thm*}{\begin{thm}}{\end{thm}}

  \newcommand{\elsonly}{}
  \newcommand{\proofbox}{\hfill $\Box$}
}

\newcommand{\baseformat}{
  \documentclass{article}

  \newenvironment{frontmatter}{}{}
  \newenvironment{keyword}{}

  \newcommand{\address}[1]{}
  \newcommand{\ead}[1]{}

  \usepackage{amsthm, amsmath}
  \theoremstyle{plain}
  \newtheorem{thm}{Theorem}
  \newtheorem{thm*}[thm]{Theorem}
  \newtheorem{conj}[thm]{Conjecture}
  \newtheorem{lem}[thm]{Lemma}
  
  \newtheorem{prop}[thm]{Proposition}
  
  \newtheorem{axiom}[thm]{Axiom}

  \theoremstyle{remark}
  \newtheorem{exmp}[thm]{Example}

  \newcommand{\bstyle}{plain}
  \newcommand{\elsonly}[1]{}
  \newcommand{\proofbox}{}
}

\baseformat

\usepackage{amsfonts,amsmath,verbatim, amscd}

\theoremstyle{plain}

\DeclareMathOperator{\Coin}{Coin}
\DeclareMathOperator{\sign}{sign}
\DeclareMathOperator{\ind}{ind}
\DeclareMathOperator{\Ind}{Ind}
\DeclareMathOperator{\id}{id}
\DeclareMathOperator{\tr}{tr}
\DeclareMathOperator{\graph}{graph}

\newcommand{\R}{\mathbb{R}}
\newcommand{\Gl}{\mathit{Gl}}
\newcommand{\adm}{\mathcal C}
\newcommand{\lto}[1]{\overset{#1}{\longrightarrow}}

\begin{document}
\bibliographystyle{\bstyle}

\begin{frontmatter}
\title{On the uniqueness of the coincidence index on orientable differentiable manifolds}
\author{P. Christopher Staecker}
\address{Messiah College, Grantham, PA.}
\ead{cstaecker@messiah.edu}

\elsonly{
}

\maketitle

\begin{abstract}
The fixed point index of topological fixed point theory is a
  well studied integer-valued algebraic invariant of a mapping which
  can be characterized by a small set of axioms. The coincidence index
  is an extension of the concept to topological (Nielsen) coincidence
  theory. We demonstrate that
  three natural axioms are sufficient to characterize the coincidence
  index in the setting of continuous mappings on oriented differentiable
  manifolds, the most common setting for Nielsen coincidence theory.
\end{abstract}
\end{frontmatter}

\section{Introduction}
For two mappings $f,g: X \to Y$, we say that $x\in X$ is a \emph{coincidence
point} of $f$ and $g$ (we write $x \in \Coin(f,g)$) if $f(x) =
g(x)$. This notion generalizes the common concept of a fixed point, in
which $X=Y$ and $g=\id$, the identity mapping. For any open set
$U\subset X$, let $\Coin(f,g,U) = \Coin(f,g) \cap U$.

The well known fixed point index of a mapping over a given subset (see
\cite{jian83} or \cite{brow71}) was
generalized to Nielsen coincidence theory  by Schirmer
\cite{schi55} in the setting of mappings between oriented manifolds
of the same (finite) 
dimension. Although Nielsen coincidence theory of nonmanifolds
\cite{gonc99}
and possibly nonorientable manifolds of nonequal dimension has been
studied, very little is known, and the coincidence index is in general
undefined in those settings.

Since definitions of the fixed point index of a mapping on a set can
be cumbersome, it is often presented in the expository literature in
terms of several properties: it is homotopy invariant, is additive on
disjoint subsets, and the total index over the whole space
equals the Lefschetz number of the mapping. This presentation is bolstered
by uniqueness results which use such properties as axioms,
and demonstrate that at most one function can satisfy the axioms (see
\cite{brow70}, in which five axioms are used). 

In \cite{fps04}, Furi,
Pera, and Spadini show that three axioms: homotopy invariance,
additivity, and a normalization axiom, are sufficient to uniquely
characterize the fixed point index in the setting of (possibly nonoriented) differentiable
manifolds. We prove a similar result in coincidence theory, though we
additionally will require our manifolds to be orientable, and we
require a stronger normalization axiom.

The approach, following the structure of 
\cite{fps04}, is to first demonstrate the uniqueness of the
coincidence index for mappings $f,g:\R^n \to \R^n$. The uniqueness
result on orientable manifolds is then obtained by approximating our mappings
first by differentiable, and then (locally) by linear maps. The last
section of the paper is devoted to the distinction between our
normalization axiom and the weaker axiom of \cite{fps04}. We show
that a uniqueness result can be obtained with the weak normalization
axiom in the special case of coincidence theory of selfmaps, but that
this uniqueness does not extend to the general setting.

I would like to thank Robert F. Brown for bringing the problem to my
attention and many helpful suggestions, Julia Weber for helpful
conversations and notes on early drafts of the paper, and the referees
for helpful suggestions which have greatly improved the quality of the paper.

\section{Preliminaries}
Throughout the paper, $X$ and $Y$ will be connected oriented differentiable
manifolds of some particular dimension $n$. For continuous maps $f,g:X
\to Y$ and some open subset $U \subset X$, we 
say that the triple $(f,g,U)$ is \emph{admissable} if the set of
coincidences of $f$ and $g$ in $U$ is compact. In particular, this
will be true provided that $U$ has compact closure and no coincidences
on its topological boundary. Let $\adm(X,Y)$ be the set of all admissable
triples of maps from $X$ to $Y$.

If $f_t, g_t: X \times [0,1] \to Y$ are homotopies and $U \subset X$
is an open subset, we say that the pair $(f_t, g_t)$ is a pair of
\emph{admissable homotopies in $U$} if
\[ \{ (x,t) \mid x \in \Coin(f_t,g_t, U), \, t \in [0,1] \} \]
is a compact subset of $X \times [0,1]$. In this case
we say that $(f_0,g_0,U)$ is \emph{admissably homotopic} to
$(f_1,g_1,U)$.

Let $D_X: H^{n-q}(X) \to H_q(X)$ and $D_Y: H^{n-q}(Y) \to H_q(Y)$ be
the Poincar\'e duality isomorphisms. Given maps $f,g: X \to Y$,
consider the composition
\[ H_q(X) \lto{f_q} H_q(Y) \lto{D_Y^{-1}} H^{n-q}(Y) \lto{g^{n-q}}
H^{n-q}(X) \lto {D_X} H_q(X), \]
where $f_q$ is the map induced in dimension $q$ homology (with
rational coefficients) by $f$, and $g^{n-q}$ is the map induced in
dimension $n-q$ cohomology by $g$. Then the \emph{Lefschetz number} of
$f$ and $g$ is defined as
\[ L(f,g) = \sum_{q = 0}^n (-1)^q \tr(D_X \circ g^{n-q} \circ D_Y^{-1}
\circ f_q). \]

Our main result concerns the uniqueness of the coincidence index with
respect to the following axioms. Throughout, $\ind$ denotes a function
$\ind:\adm(X,Y) \to \R$.

\begin{axiom}[Additivity axiom]
Given $(f,g,U)\in \adm(X,Y)$, if $U_1$ and $U_2$ are disjoint open
subsets of $U$ with $\Coin(f,g, U) \subset U_1 \cup U_2$, then
\[ \ind(f,g,U) = \ind(f,g,U_1) + \ind(f,g,U_2). \]
\end{axiom}

\begin{axiom}[Homotopy axiom]
If $(f_0, g_0, U) \in \adm(X,Y)$ and $(f_1, g_1, U)\in \adm(X,Y)$ are
admissably homotopic, then
\[ \ind(f_0, g_0, U) = \ind(f_1, g_1, U). \]
\end{axiom}

\begin{axiom}[Normalization axiom]
For $(f,g,X) \in \adm(X,Y)$, we have
\[ \ind(f,g,X) = L(f,g). \]
\end{axiom}

\begin{axiom}[Weak normalization axiom]
If $(f,\id,X) \in \adm(X,X)$, with $f$ a constant map and $\id$ the
identity map, then $\ind(f,\id,X) = 1$.
\end{axiom}

The weak normalization axiom is a special case of the 
normalization axiom where $f$ and $g$ are selfmaps, since a
straightforward calculation will show that $L(f,\id) = 1$ when $f$ is
constant. The full normalization axiom is needed to handle coincidence
theory of non-selfmaps, in which case the weak normalization axiom
will not apply. Many of our results, however, require only the weak
version, and we will be careful to indicate in our hypotheses which
normalization axiom is being used.

Our main result is the following:
\begin{thm} \label{mainthm}
There is at most one function $\ind:\adm(X,Y) \to \R$ satisfying the
additivity, homotopy, and normalization axioms.
\end{thm}

This result is a slightly weakened generalization of the main result
of \cite{fps04}, which demonstrates the uniqueness of the fixed point
axiom with respect to the additivity, homotopy, and weak normalization
axioms. 

We can derive some immediate corollaries from these three axioms,
which will be useful in the exposition to follow. In
the following propositions, $\ind: \adm(X,Y) \to \R$ denotes any
function which satisfies the homotopy axiom, additivity axiom, and,
whenever $X=Y$, the weak normalization axiom.

\begin{prop}[Fixed-point property] \label{fpprop} If $(f,\id,U) \in
\adm(X, X)$, then
\[ \ind(f,\id,U) = i(f,U), \]
where the right-hand side is the well known fixed point index. 
\end{prop}
\begin{proof}
To demonstrate that $\ind(f,\id,U)$ is the fixed point index, we need
only show that it obeys the three axioms described in
\cite{fps04}. Each of those three axioms are clearly special cases of
our homotopy, additivity, and weak normalization axioms, so they are
satisfied.
\proofbox
\end{proof}

\begin{prop} [Empty set property] The triple $(f,g,\emptyset)$ is admissable,
and has index 0.
\end{prop}
\begin{proof}
By the additivity axiom, we write $\emptyset= \emptyset \cup
\emptyset$, and we have
\[ \ind(f,g,\emptyset) = \ind(f,g,\emptyset) + \ind(f,g,\emptyset), \]
and so $\ind(f,g,\emptyset) = 0$.
\proofbox \end{proof}

\begin{prop}[Excision property] If $(f,g,U) \in \adm(X,Y)$ and $U'\subset U$
is an open subset containing $\Coin(f,g,U)$, then
$\ind(f,g,U) = \ind(f,g,U')$.
\end{prop}
\begin{proof}
By the additivity axiom, we have:
\[ \ind(f,g,U) = \ind(f,g,U') + \ind(f,g, \emptyset) = \ind(f,g,U'). \]
\proofbox \end{proof}

\begin{prop}[Solution property] If $\ind(f,g,U) \neq 0$, then $f$ and $g$
have a coincidence on $U$.
\end{prop}
\begin{proof}
We prove the contrapositive: if $f$ and $g$ have no coincidence on
$U$, then by the excision property we have
\[ \ind(f,g,U) = \ind(f,g,\emptyset) = 0. \]
\proofbox \end{proof}

\section{The coincidence index for mappings of $\R^n$}
Let $M_n$ be the space of linear maps $A:\R^n \to
\R^n$, and $\Gl_n$ the subspace of invertible maps. Define the set $N
\subset M_n \times M_n$ as the set of all pairs of maps $(A,B)$ such
that $B-A \in \Gl_n$. Note that $N$ gives precisely the linear mappings for
which $(A,B,\R^n)\in \adm(\R^n, \R^n)$, since $\Coin(A,B)$ will be a linear
subspace of $\R^n$, and thus will be compact if and only if it is
$\{0\}$, in which case $\det(B-A) \neq 0$.

We define 
\begin{align*}
N^+ &= \{(A,B) | \det(B-A) > 0\}, \\
N^- &= \{(A,B) | \det(B-A) < 0\},
\end{align*}
and note that $N$ is a disconnected set, with components $N^+$ and
$N^-$. Each of these is open in $M_n\times M_n$, and
therefore path-connected. Thus if some function $\ind:\adm(\R^n,
\R^n) \to \R$ satisfies the homotopy axiom, then for any $(A,B)
\in M_n \times M_n$, the value of $\ind(A,B,\R^n)$
depends only on the path component containing $(A,B)$.

If $\ind$ additionally satisfies the weak normalization axiom, then
$\ind(0,\id,\R^n) = 1$, and so we have that $\ind(A,B,\R^n) = 1$ for
all $(A,B) \in N^+$. 

\begin{lem}
For any $(A,B) \in N^-$, if $\ind:\adm(\R^n,\R^n) \to \R$ satisfies
the additivity, homotopy, and weak normalization axioms, we have
\[ \ind(A,B,\R^n) = -1. \]
\end{lem}
\begin{proof}
It suffices to show that $\ind(A,B,\R^n) = -1$ for some particular
pair $(\hat A,\hat B) \in N^-$. We will show this for $\hat B = \id$ and
$\hat A$ the linear map of the same name described in Lemma 3.1 of
\cite{fps04}: 
\[ \hat A(x_1, \dots, x_{n-1}, x_n) = (0, \dots, 0, 2x_n). \]
Since we are taking $\hat B$ to be the identity map, Proposition
\ref{fpprop} implies that we need only demonstrate that the fixed
point index of $\hat A$ is $-1$, and this is demonstrated in Lemma 3.1
of \cite{fps04}. 
\proofbox \end{proof}

The above discussion and lemma give our first uniqueness result.
\begin{lem}\label{linearindex}
For any $\ind:\adm(\R^n, \R^n)\to \R$ satisfying the additivity,
homotopy, and weak normalization axioms, if $A,B: \R^n \to \R^n$ are linear
maps and $(A,B,\R^n) \in \adm(\R^n, \R^n)$, then
\[ \ind(A,B,\R^n) = \sign \det(B-A). \]
\end{lem}

For differentiable maps $f,g: X \to Y$, let $df_x, dg_x: \R^n \to
\R^n$ denote the 
derivative maps of $f$ and $g$ at the point $x\in X$. If the triple
$(f,g,U)$ is admissable, and $dg_p - df_p \in \Gl_n$ for every $p \in
\Coin(f,g,U)$, then we say that the triple $(f,g,U)$ is
\emph{nondegenerate}.  

Since we have established that the index is very well behaved 
for linear mappings, the following linearization result for the index
is very useful.

\begin{lem}\label{isolation}
Let $(f,g,U)\in \adm(\R^n, \R^n)$ be a nondegenerate
triple. Then each coincidence point
$p$ of $f$ and $g$ is 
isolated, and for some isolating neighborhood $V$ of $p$, we have
\[ \ind(f,g,V) = \ind(df_p, dg_p, \R^n), \]
where $\ind:\adm(\R^n, \R^n) \to \R$ is any function satisfying the
additivity, homotopy, and weak normalization axioms.
\end{lem}
\begin{proof}
Let $p$ be a coincidence point. Since $f$ and $g$ are differentiable
at $p$, we have 
\begin{align*}
f(x) &= f(p) + df_p(x-p) + |x-p|\epsilon(x-p), \\
g(x) &= g(p) + dg_p(x-p) + |x-p|\delta(x-p),
\end{align*}
where $\epsilon, \delta: U \to \R^n$ are continuous with $\epsilon(0)
= \delta(0) = 0$. Then we have
\begin{align*}
|g(x) - f(x)| &= \left|(dg_p - df_p)(x-p) + |x-p|(\delta(x-p) -
 \epsilon(x-p))\right| \\
&\geq |(dg_p - df_p)(x-p)| - |x-p||\delta(x-p) - \epsilon(x-p)|
\end{align*}
Since $dg_p - df_p$ is a linear map, we have
\begin{align*}
|(dg_p - df_p)(x-p)| &= |x-p|\left|(dg_p -
df_p)\left(\frac{x-p}{|x-p|}\right)\right| \\
&\geq  |x-p|\left(\inf_{|v|=1} |(dg_p - df_p)(v)|\right),
\end{align*}
and so
\[ |g(x) - f(x)| \geq |x-p|\left( \inf_{|v|=1}|(dg_p - df_p)(v)| -
|\delta(x-p) - \epsilon(x-p)| \right). \]
The infimum is known to be strictly positive, since $(f,g,U)$ is
nondegenerate. Thus for all $x$ sufficiently close to $p$, the term
$|\delta(x-p) - \epsilon(x-p)|$ will be sufficiently small so that
$|g(x) - f(x)|>0$, which is to say that $p$ is an isolated coincidence
point. 

It remains to show that if $V$ is an isolating neighborhood of $p$, then 
\[ \ind(f,g,V) = \ind(df_p, dg_p, \R^n). \]
Define the following homotopies:
\begin{align*}
f_t(x) &= f(p) + df_p(x-p) + t|x-p|\epsilon(x-p), \\
g_t(x) &= g(p) + dg_p(x-p) + t|x-p|\delta(x-p).
\end{align*}
By the same argument as above, for all $t\in [0,1]$ we have
\[ |g_t(x) - f_t(x)| > 0 \]
for all $x$ in some neighborhood $W \subset V$ of $p$. Thus
$(f_t,g_t)$ is an admissable homotopy on $W$ of $f$ and $g$ to the
affine linear maps $x \mapsto f(p) + df_p(x-p)$ and $x \mapsto g(p) +
dg_p(x-p)$. These affine linear maps have $p$ as their only
coincidence point, since $dg_p - df_p \in \Gl_n$. Thus we have
\begin{align*} 
\ind(f,g,W) &= \ind(f(p) + df_p(x-p), g(p) + dg_p(x-p), W) \\
&= \ind(f(p) + df_p(x-p), g(p) + dg_p(x-p), \R^n) 
\end{align*}
by the homotopy axiom and the excision property.

But these affine maps are clearly admissably homotopic to
$df_p$ and $dg_p$, and so we have
\[ \ind(f,g,V) = \ind(f,g,W) = \ind(df_p,dg_p,\R^n), \]
where the first equality is by the excision property, and the second
is by the homotopy axiom. 
\proofbox \end{proof}

\section{The coincidence index for mappings of orientable manifolds}
Lemmas \ref{linearindex} and \ref{isolation} established
the uniqueness of the index for nondegenerate pairs of mappings on
$\R^n$ with respect to the additivity, homotopy, and weak
normalization axioms. The next two lemmas show that any index for nondegenerate
pairs on arbitrary orientable manifolds is computed similarly to the
index of Lemma \ref{linearindex}. 

\begin{lem} \label{mfldisolation} If $(f,g,U)\in \adm(X,Y)$ is a nondegenerate
triple, then each coincidence point is isolated.
\end{lem}
\begin{proof}
For some $p \in \Coin(f,g,U)$, choose a chart $Z_q \subset Y$ containing
$q = f(p) = g(p)$ with a diffeomorphism
$\psi_q: Z_q \to \R^n$. Choose a chart $W_p \subset f^{-1}(Z_q) \cap
g^{-1}(Z_q)$ containing $p$ with a diffeomorphism $\phi_p:W_p \to \R^n$.
Define $\omega: \adm(W_p, X) \to \adm(\R^n, \R^n)$ by
\[ \omega (f,g,W_p) = (\psi_q \circ f \circ \phi^{-1}_p, \psi_q \circ
g \circ \phi^{-1}_p, \phi_p(W_p)). \]

Since $(f,g,U)$ is nondegenerate, then clearly $(f,g,W_p)$ will be
nondegenerate. For any  
coincidence point $x \in \Coin(\omega(f,g,W_p))$, we have
\[ d(\psi_q \circ f \circ \phi^{-1}_p)_x - d(\psi_q \circ g \circ
\phi^{-1}_p)_x = d(\psi_q)_{f(\phi^{-1}_p(x))} (df_{\phi^{-1}_p(x)} -
dg_{\phi^{-1}_p(x)}) d(\phi^{-1}_p)_x, \]
since $f(\phi^{-1}_p(x)) = g(\phi^{-1}_p(x))$.
Since $\phi_p$ and $\psi_q$ are diffeomorphisms, the above will be in
$\Gl_n$ by nondegeneracy of $(f,g,W_p)$. Thus 
$\omega(f,g,W_p)$ is nondegenerate, which implies by Lemma
\ref{isolation} that $\Coin(\omega(f,g,W_p))$ is a set of
isolated points. But this coincidence set is in diffeomorphic
correspondence via $\phi_p$ and $\psi_q$ to $\Coin(f,g,W_p)$, so this
too is a set of isolated points, and in particular $p$ is an
isolated coincidence point.
\proofbox \end{proof}

We note that the next Lemma contains our only use of the full (rather
than the weak) normalization axiom, and also our only use of the
orientability hypothesis on the manifolds $X$ and $Y$.

\begin{lem} \label{nondegindex} If $(f,g,U)\in \adm(X,Y)$ is a nondegenerate
  triple, and $\ind:\adm(X,Y) \to \R$ satisfies the additivity,
  homotopy, and normalization axioms, then
\begin{equation}\label{signofdet} \ind(f,g,U) = \sum_{p \in
	\Coin(f,g,U)} \sign(\det(dg_p - df_p)).
\end{equation}
\end{lem}
\begin{proof}
Note that since $(f,g,U)$ is nondegenerate, by Lemma \ref{mfldisolation} the
coincidence set $C = \Coin(f,g,U)$ consists of finitely many isolated
points. 

First we prove the theorem in the special case where $X=Y=\R^n$.
For each coincidence point $p\in U$,
let $V_p$ be an isolating neighborhood of $p$. Since $C$ is compact,
we may choose the sets $V_p$ to be pairwise 
disjoint. Then by the additivity property, Lemma \ref{isolation}, and Lemma
\ref{linearindex}, we have
\begin{align*} 
\ind(f,g,U) &= \sum_{p \in C} \ind(f,g,V_p) = \sum_{p \in C}
\ind(df_p,dg_p,\R^n) \\ 
&= \sum_{p \in C} \sign(\det(dg_p - df_p)). 
\end{align*}

It remains to prove the result in the general case that $X$ and $Y$
are not both $\R^n$.
About each isolated coincidence point $p$, choose
pairwise disjoint isolating neighborhoods $W_p$ diffeomorphic to $\R^n$
by orientation preserving diffeomorphisms $\phi_p: W_p \to
\R^n$. We also choose neighborhoods $Z_q$ of each coincidence value $q
= f(p) = g(p)$ with $Z_q$ diffeomorphic to $\R^n$ by orientation
preserving diffeomorphisms $\psi_q: Z_q \to \R^n$.
Let $\omega:\adm(W_p, Y) \to \adm(\R^n, \R^n)$ be as in Lemma
\ref{mfldisolation}. 

Note that $\omega$ has an inverse given by 
\[ \omega^{-1} (F,G,S) = (\psi^{-1}_q \circ F \circ \phi_p,
\psi_q^{-1} \circ G \circ \phi_p, \phi^{-1}_p(S)), \]
where $(F,G,S) \in \adm(\R^n, \R^n)$. Then we have trivially that 
\[ \ind(f,g,W_p) = \ind \circ \omega^{-1} (\omega(f,g,W_p)). \]

We now note that $\ind \circ \omega^{-1}: \adm(\R^n, \R^n) \to \R$
satisfies the additivity, homotopy, and normalization axioms. The
additivity and homotopy axioms are 
clear, but a brief calculation is needed for the normalization
axiom. We must show that $\ind \circ \omega^{-1}(F, G, \R^n) =
L(F,G)$.

We have
\begin{align*}
\ind(\omega^{-1}(F,G,\R^n)) &= \ind(\psi_q^{-1} \circ F \circ
\phi_p, \psi_q^{-1} \circ G \circ \phi_p, W_p) \\
&= L(\psi_q^{-1} \circ F \circ
\phi_p, \psi_q^{-1} \circ G \circ \phi_p) 
\end{align*}
by the normalization axiom. Consider the diagram:
\[ \begin{CD}
H_*(W_p) @>\phi_{p*}>> H_*(\R^n) @>F_*>> H_*(\R^n) @>\psi_{q*}^{-1}>>
H_*(Z_q) \\
@AD_X^{-1}AA  @AD^{-1}_{\R^n}AA @VD_{\R^n}VV @VD_YVV \\
H^{n-*}(W_p) @<\phi_p^*<< H^{n-*}(\R^n) @<G^{n-*}<< H^{n-*}(\R^n)
@<\psi_q^{-1*}<< H^{n-*}(Z_q)
\end{CD} \]
where starred maps are the induced maps in (co)homology. The Lefschetz
number above is the alternating sum of the traces of the maps formed
by a clockwise walk around the perimeter of the diagram, starting at
the upper left corner.

Note by the functoriality of the duality maps that the squares at
right and left in the diagram will commute. Using this fact, along
with cyclic permutation of maps inside the trace, we have
\begin{align*}
\tr(D_X^{-1} \circ \phi_p^* \circ \psi_q^{-1*} \circ D_Y \circ
\psi_{q*}^{-1} \circ F_* \circ \phi_{p*}) &= \tr(D_X^{-1} \circ
\phi_p^* \circ G^* \circ D_{\R^n} \circ F_* \circ \phi_{p*}) \\
&= \tr(\phi_{p*} \circ D_X^{-1} \circ \phi_p^* \circ G^* \circ
D_{\R^n} \circ F_*) \\
&= \tr(D_{\R^n}^{-1} \circ G^* \circ D_{\R^n} \circ F_*), 
\end{align*}
and so 
\[ L(\psi_q^{-1} \circ F \circ
\phi_p, \psi_q^{-1} \circ G \circ \phi_p) = L(F,G), \]
which shows that $\ind \circ \omega^{-1}$ satisfies the normalization
axiom.

We have shown that $\ind \circ \omega^{-1}$ is a real valued function
on $\adm(\R^n, \R^n)$ which satisfies the three axioms. Thus, by the
special case above, it is calculated according to (\ref{signofdet}).

Thus we have:
\begin{align*}
\ind(f,g,U) &= \sum_{p \in C} \ind(f,g,W_p) = \sum_{p\in C}
\ind\circ \omega^{-1}(\omega(f,g,W_p)) \\
&= \sum_{p\in C} \sign (\det(d\psi_q \cdot dg_p \cdot
d\phi_p^{-1} - d\psi_q \cdot df_p \cdot d\phi_p^{-1}))
\end{align*} 
by the derivative chain rule. (For brevity, we have written
$d\phi_p^{-1}$ to indicate $d(\phi_p^{-1})_p$ and 
$d\psi_q$ to indicate $d(\psi_q)_q$.) But the above is simply
\[ \sum_{p\in C} \sign(\det( d\psi_q (dg_p - df_p)
d\phi_p^{-1})) = \sum_{p\in C} \sign(\det(dg_p - df_p)) \]
since all $\phi_p$ and $\psi_q$ are taken to be orientation preserving.
\proofbox \end{proof}

One final Lemma is required, showing that any pair of maps can be
approximated by a nondegenerate pair of maps.

\begin{lem} \label{nondegapprox} Let $(f,g,U)\in \adm(X,Y)$, and $V
\subset U$ be an open subset containing $\Coin(f,g, U)$ with compact
closure $\bar V \subset U$. Then $(f,g,V)$ is admissably
homotopic to a nondegenerate triple.  
\end{lem}
\begin{proof}
In order to facilitate the construction of explicit straight-line
homotopies, we begin by embedding our manifolds in some Euclidean
space, and approximating our maps accordingly by close polynomial
approximations. 

Without loss of generality, assume that $X$ and $Y$ are embedded in $\R^k$ for
some $k$. By the $\epsilon$-neighborhood theorem (see \cite{gp74}),
there exists a neighborhood $\Omega$ of $Y$ in $\R^k$ with a
submersive retraction $r:\Omega \to Y$ such that $|x-r(x)| =
\mathrm{dist}(x,Y)$.

For any $\delta > 0$, by the Weierstrass approximation theorem, there
exist polynomial maps $f_\delta,g_\delta:\R^k \to \R^k$ with
$|f(x)-f_\delta(x)| < \delta$ and $|g(x)-g_\delta(x)|<\delta$ for all
$x \in \bar V$. Since $\bar V$ is compact, we may choose $\delta$ 
sufficiently small so that the homotopies
\begin{align*} 
f_t(x) &= r((1-t)f(x) + tf_\delta(x)) \\
g_t(x) &= r((1-t)g(x) + tg_\delta(x))
\end{align*}
are well-defined and $(f_t,g_t)$ is an admissable homotopy on
$V$. Thus the pair $(f,g)$ is admissably homotopic to
$(r\circ f_\delta, r \circ g_\delta)$ in $V$. To simplify our
notation, let $f' = r \circ f_\delta$, and $g' = r \circ g_\delta$.
We will show that the pair $(f',g')$ is
admissably homotopic to some nondegenerate pair.

Let $B$ be a ball about the origin sufficiently small that for all
$x\in \bar V$ and $y \in B$ the function $x \mapsto r(f'(x) + y)$ is
defined and has no coincidences with $g'$ on $\partial V$. Then define
$H:V \times B \to V \times Y$ by
\[ H(x,y) = (x, r(f'(x) + y)), \]
and note that the derivative map
\[ dH_{(x,y)}: \R^k \times \R^k \to \R^k \times T_{r(f'(x) + y)} Y \]
is surjective (since $r$ is a submersion). Thus at points where the
image of $H$ intersects with the graph of $g'$ in $V \times Y$, the
intersection will be transversal.

By the transversality theorem (see \cite{gp74}), $H(x,y)$ is
transversal to $\graph g'$ for almost all $y\in B$. Choose one such
$\bar y$ so that
\[ \bar H(x) = (x, r(f'(x) + \bar y)) \]
is transversal to $\graph g'$. Note that a pair of mappings is
nondegenerate if and only if their graphs intersect transversally.
This means that $(r(f' + \bar y), g', \bar V))$ is a nondegenerate triple.

But our assumption on $B$ means that the homotopy given by
\[ (r(f' + t\bar y), g') \]
is admissable on $\bar V$, and thus the nondegenerate triple above is
admissably homotopic to $(f', g', \bar V)$,
which has already been shown to be admissably homotopic to $(f,g,V)$.
\proofbox \end{proof}

Now we are ready to prove Theorem \ref{mainthm}:
\begin{proof}
If $(f,g,U) \in \adm(X,Y)$, then by Lemma \ref{nondegapprox} there is an
open set $V \subset U$ containing $\Coin(f,g,U)$ with 
$(f,g,V)$ admissably homotopic to a nondegenerate triple $(f',g',V)$. Thus
by excision and the homotopy axiom, we have $\ind(f,g,U) =
\ind(f',g',U)$, and so by Lemma \ref{nondegindex} we have
\[ \ind(f,g,U) = \sum_{p \in \Coin(f',g', V)} \sign(\det(dg'_p -
df'_p)). \]

The above calculation does not depend on the choice of
nondegenerate triple $(f',g',U)$, since any alternative choice would
automatically have the same index by the homotopy axiom. Since any
coincidence index must obey the calculation given above, we have shown
that all indices must agree, and that there is at most one.
\proofbox \end{proof}

At this point we will briefly address the issue of the existence of
the coincidence index. The above lemmas could be used to define a
coincidence index in the following way: given any admissable triple
$(f,g,U)$, let $(f',g',V)$ be a nondegenerate triple given by Lemma
\ref{nondegapprox}, admissably homotopic to $(f,g,V)$ with
$\Coin(f,g,U) \subset V$. We then define $\ind(f,g,U) =
\ind(f',g',V)$, where $\ind(f',g',V)$ is computed according to
(\ref{signofdet}).

An index defined in this way would clearly satisfy our axioms,
provided that it is well defined. We believe this to be the case, but
a verification is apparently nontrivial, and would require a
demonstration of: 
\begin{conj}
If $(f,g,U)$ and $(h,k,U)$ are admissably homotopic
nondegenerate triples, then
\[ \sum_{p\in \Coin(f,g,U)} \sign(\det(dg_p - df_p)) = \sum_{q \in
\Coin(h,k,U)} \sign(\det(dk_q - dh_q)). \]
\end{conj}

\section{Uniqueness and the weak normalization axiom}
It is natural to ask whether any meaningful version of the above
theorem can be proved using the weak normalization axiom in place of
the full normalization axiom. Since the weak normalization can only
apply to selfmaps, an obvious question is: is there is a unique
coincidence index on $\adm(X, X) \to \R$ which satisfies the homotopy,
additivity, and weak normalization axioms? We answer this question in
the affirmative, and additionaly show that orientability of $X$ is not
required in this setting.

\begin{thm} \label{selfmapthm}
For any particular (perhaps nonorientable) differentiable manifold $X$,
there is at most one coincidence index $\ind: \adm(X,X) \to \R$
satisfying the additivity, homotopy, and weak normalization axioms.
\end{thm}
\begin{proof}
Note that the proof of Lemma \ref{nondegindex} is the only place where
either of the orientability hypothesis or the full (as opposed to the weak)
normalization axiom is used. Thus our proof here may make use of any
of our lemmas except for Lemma \ref{nondegindex}.

As in the proof of Theorem \ref{mainthm}, our proof consists of
deriving a formula 
for the computation of the index of some triple $(f,g,U)$. By Lemma
\ref{nondegapprox} we may assume without loss of generality that
$(f,g,U)$ is nondegenerate. As in the proof of Lemma
\ref{nondegindex}, an explicit formula for
the index is clear in the special case where $X= \R^n$. In the case
where $X$ is not $\R^n$, we first change the triple $(f,g,U)$ by a
homotopy, and then follow similar steps to those used in Lemma
\ref{nondegindex}. 

For some isolated coincidence point $p$, let $q = f(p) = g(p)$. Let
$\gamma:[0,1] \to X$ be a path in $X$ from $q$ to $p$ which avoids all other
points of $\Coin(f,g,U)$, and let $V$ be a contractible neighborhood of
$\gamma$ with closure disjoint from $\Coin(f,g,U)$. Since $V$ is
homeomorphic to an open ball in $\R^n$, there are homotopies $f_t, g_t$ 
such that $f_t$ and $g_t$ agree for all $t$ with $f$ and $g$, respectively, on
$U-V$, and $\Coin(f_t,g_t,V) = \{p\}$ with $f_t(p) = g_t(p) =
\gamma(t)$. Thus the triple $(f_t,g_t,U)$ will be an admissable
homotopy of the triple $(f,g,U)$ to some triple $(f_1,g_1,U)$, with
$f_1(p) = g_1(p) = p$. 

In this way, we have converted the coincidence point at $p$ into a
fixed point at $p$ without disturbing the behavior of $f$ and $g$ at
the other coincidence points. By iterating the above construction for
each coincidence point $p$, we obtain an admissable triple $(f', g',
U)$ admissably homotopic to $(f,g,U)$ such that all points of
$\Coin(f', g', U)$ are actually fixed points of both $f'$ and
$g'$.

We now continue as in the proof of Lemma \ref{nondegindex}. About each
isolated coincidence point $p \in \Coin(f',g',U)$, choose pairwise
disjoint isolating neighborhoods $W_p$ diffeomorphic to $\R^n$ by
diffeomorphisms $\phi_p: W_p \to \R^n$. Letting $V_p = f'^{-1}(W_p)
\cap g'^{-1}(W_p) \cap W_p$, we will define $\omega(f',g',V_p) \in
\adm(\R^n, \R^n)$ as in Lemma \ref{mfldisolation},
only this time our two charts can be taken to be the same, since $p$
is a fixed point. Let  
\[ \omega(f',g',V_p) = (\phi_p \circ f' \circ \phi^{-1}_p, \phi_p
\circ g' \circ \phi_p^{-1}, \phi(V_p)), \]
and as in Lemma \ref{nondegindex}, we note that $\omega$ has an
inverse given by
\[ \omega^{-1}(F,G,S) = (\phi_p^{-1} \circ F \circ \phi_p, \phi_p^{-1}
\circ G \circ \phi_p, V_p), \]
and trivially we have $\ind(f',g',V_p) = \ind \circ \omega^{-1}
(\omega(f',g',V_p))$. 

Now we note that $\ind \circ \omega^{-1}:
\adm(\R^n, \R^n) \to \R$ satisfies the additivity, homotopy, and weak
normalization axioms. The additivity and homotopy axioms are clear,
and the weak normalization is not difficult: if $F$ is the
constant function $F(x) = c$ and $G$ is the identity, then
\[ \ind \circ \omega^{-1} (F, G, S) = \ind(\phi^{-1}(c), \id, V_p), \]
and this equals 1 by the weak normalization axiom.

Thus $\ind \circ \omega^{-1}$ satisfies the additivity, homotopy, and
weak normalization axioms, and so is the coincidence index by our
special case above for $X = \R^n$. Thus, letting $C = \Coin(f',g',U)$,
we have
\begin{align*} 
\ind(f,g,U) &= \ind(f',g',U) = \sum_{p\in C}
\ind(f',g',V_p) = \sum_{p \in C} \ind \circ \omega^{-1}(\omega(f',g',V_p)) \\
&= \sum_{p \in C} \ind(\omega(f',g',V_p)) = \sum_{p \in C}
\sign(\det(d\phi_p \cdot dg'_p \cdot d\phi_p^{-1} - 
d\phi_p \cdot df'_p \cdot d\phi_p^{-1})) \\
&= \sum_{p \in C} \sign(\det(dg'_p - df'_p))
\end{align*}
By the homotopy axiom, this formula is independent of the choice of
the admissible homotopy to $(f',g',U)$, and the uniqueness is shown.
\end{proof}

The above theorem could be restated as follows: If $\adm_s$ is the
class of \emph{all} admissable triples of selfmaps on differentiable
manifolds, then there is a unique coincidence index $\ind: \adm_s \to
\R$ which satisfies the additivity, homotopy, and weak normalization
axioms. This index, for any triple $(f,g,U)$, would be defined by using
the unique index defined on $\adm(X,X)$, where $X$ is the domain of
$f$ and $g$. Our above theorem, in this sense, can be seen as a direct
generalization of the result of \cite{fps04} to coincidence theory of
selfmaps. 

A natural question to ask is whether a further extension of
\cite{fps04} can be made to non-selfmaps
as follows: If $\adm$ is the set of \emph{all} admissable triples
$(f,g,U)$ where $f$ and $g$ are maps between orientable differentiable
manifolds of the same dimension, is a coincidence
index satisfying the additivity and homotopy axioms, which additionally
satisfies the weak normalization axiom whenever $f$ and $g$ are
selfmaps, unique? Such a uniqueness result would be stronger than our Theorem
\ref{mainthm} but is false as the following example illustrates.

\begin{exmp}
Let $\adm$ be as above, the set of all admissable triples. Letting
$\Ind_{(X,Y)}$ denote the unique coincidence index on $\adm(X,Y)$
given by Theorem \ref{mainthm}, we can define a single index
$\Ind:\adm \to \R$ by $\Ind(f,g,U) = \Ind_{(X,Y)}(f,g,U)$
when $(f,g,U) \in \adm(X,Y)$.
This gives a single coincidence index $\Ind$
on $\adm$ which satisfies the additivity, homotopy, and normalization
axioms. 

For any $c \in \R$, define $i_c:\adm \to \R$ as follows:
\[ i_c(f,g,U) = \begin{cases} \Ind(f,g,U) \quad & \text{if $f$, $g$
are selfmaps,} \\
c \cdot \Ind(f,g,U) & \text{otherwise}. \end{cases} \]
For any value of $c$, this function clearly satisfies the additivity
and homotopy axioms, and also satisfies the weak normalization axiom
in the case when $f$ and $g$ are selfmaps.
\end{exmp}

Variation of the parameter $c$ in the above example will produce many
distinct ``indices'' on $\adm$ which satisfy the additivity, homotopy,
and (when applicable) weak normalization axioms. Thus Theorems
\ref{mainthm} and \ref{selfmapthm} seem to be the best uniqueness
results available in our setting.

\end{document}